\documentclass{amsart}

\usepackage[normalem]{ulem}

\usepackage{lineno}

\usepackage{caption}

\usepackage{soul}

\usepackage{graphicx}

\usepackage{amssymb, amsmath, amsthm, hyperref, euscript, color}
\usepackage[dvipsnames]{xcolor}

\usepackage{amscd}

\usepackage{bbm}

\usepackage{enumitem}

\usepackage[english]{babel}

\numberwithin{equation}{section}

\theoremstyle{plain}

\newtheorem{theorem}{Theorem}

\newtheorem{corollary}[theorem]{Corollary}

\newtheorem{proposition}[theorem]{Proposition}

\newtheorem{remark}[theorem]{Remark}

\newtheorem{lemma}[theorem]{Lemma}

\newtheorem{thm}{Theorem}

\theoremstyle{definition}

\newtheorem{example}[theorem]{Example}

\newcommand{\R}{\mathbb{R}}

\newcommand{\Z}{\mathbb{Z}}

\newcommand{\N}{\mathbb{N}}

\newcommand{\id}{\textnormal{id}}

\newcommand{\rot}{\textnormal{rot}}

\DeclareMathOperator{\Gl}{GL}

\DeclareMathOperator{\Pin}{Pin}

\makeatletter
\@namedef{subjclassname@2020}{%
  \textup{2020} Mathematics Subject Classification}
\makeatother

\begin{document}

\author{Valentina Bais and Rafael Torres}

\title[Smooth structures on 4-manifolds via twisting operations]{Smooth structures on non-orientable 4-manifolds via twisting operations}

\address{Scuola Internazionale Superiore di Studi Avanzati (SISSA)\\ Via Bonomea 265\\34136\\Trieste\\Italy}

\email{$\{vbais, rtorres\}$@sissa.it}



\maketitle

\emph{Abstract}: Five observations compose the main results of this note. The first records the existence of a smoothly embedded 2-sphere $S$ inside $\R P^2\times S^2$ such that performing a Gluck twist on $S$ produces a manifold $Y$ that is homeomorphic but not diffeomorphic to the total space of the non-trivial 2-sphere bundle over the real projective plane $S(2\gamma \oplus \R)$. The second observation is that there is a 5-dimensional cobordism with a single 2-handle between the 4-manifold $Y$ and a mapping torus that was used by Cappell-Shaneson to construct an exotic $\R P^4$. This construction of $Y$ is similar to the one of the Cappell-Shaneson homotopy 4-spheres. The third observation is that twisting an embedded real projective plane inside $Y$ produces a manifold that is homeomorphic but not diffeomorphic to the circle sum of two copies of $\R P^4$. The fourth observation records new examples of pairs of homeomorphic but not diffeomorphic closed 4-manifolds with Euler characteristic one. These include the total space of the nontrivial $\R P^2$-bundle over $\R P^2$. Knotting phenomena of 2-spheres in non-orientable 4-manifolds that stands in glaring contrast with known phenomena in the orientable domain is pointed out in the fifth observation.

\section{Introduction}\label{Introduction}

Akbulut \cite{[AkbulutA], [Akbulut2]}, \cite[Section 9.5]{[Akbulut3]} showed the existence of 2-spheres\begin{center}$S\subset (S^3\widetilde{\times} S^1)\# (S^2\times S^2)$ and $S'\subset (S^3\widetilde{\times} S^1)\# \mathbb{CP}^2\#\overline{\mathbb{CP}^2}$\end{center} such that performing a Gluck twist to either of them yields a 4-manifold that is homeomorphic but not diffeomorphic to the connected sum of the non-trivial 3-sphere bundle over the circle $S^3\widetilde{\times} S^1$ with a copy of $S^2\times S^2$. The first main result of this note records such examples in 4-manifolds that are known to be irreducible. Recall that a smooth 4-manifold $X$ is said to be irreducible if for every smooth connected sum decomposition $X = X_1\# X_2$, either $X_1$ or $X_2$ is homeomorphic to the 4-sphere \cite[Definition 10.1.17]{[GompfStipsicz]}. 

\begin{thm}\label{Theorem A}There is a smoothly embedded 2-sphere $S\subset \R P^2\times S^2$ with trivial normal bundle $\nu( S) = D^2\times S^2$ such that the complements\begin{center}$(\R P^2\times S^2)\setminus \nu(\{pt\}\times S^2)$ and $(\R P^2\times S^2)\setminus \nu (S)$\end{center}are not homotopy equivalent, and doing a Gluck twist on $S$ yields a 4-manifold $Y$ that is homeomorphic but not diffeomorphic to the non-orientable total space $S(2\gamma \oplus \R)$ of the non-trivial 2-sphere bundle over the real projective plane.

\end{thm}

Performing a Gluck twist on the 2-sphere $\{pt\}\times S^2\subset \R P^2\times S^2$ yields $S(2\gamma \oplus \R)$ as explained in Section \ref{Section Gluck Twists}. The notation $S(2\gamma \oplus \R)$ is explained in Section \ref{Section Bundles}. A key ingredient in our proof of Theorem \ref{Theorem A} is a result due to Akbulut \cite[Section 4]{[Akbulut1]}, which we state in Theorem \ref{Theorem Akbulut}. The inequivalent smooth structure $Y$ of Theorem \ref{Theorem A} was constructed in \cite[Theorem 1]{[Torres]}; cf. Proposition \ref{Proposition Sphere Bundle}. Our next result provides a different perspective on the production of the 4-manifold $Y$ of Theorem \ref{Theorem A}. 

\begin{thm}\label{Theorem B} Let $M_A$ be the Cappell-Shaneson's non-orientable 3-torus bundle over $S^1$ with monodromy $A$ as in (\ref{Matrix 1}), and let $\alpha_0\subset M_A$ be a cross-section. Let $\alpha_0^2\subset M_A$ be an orientation-preserving simple loop that represents the element $[\alpha_0]^2\in \pi_1(M_A)$. 

$\bullet$ Surgery on $\alpha_0^2$ with the canonical framing yields $\R P^2\times S^2$.

$\bullet$ Surgery on $\alpha_0^2\subset M_A$ with the non-canonical framing yields the 4-manifold $Y$ of Theorem \ref{Theorem A}, which is homeomorphic but not diffeomorphic to $S(2\gamma \oplus \R)$.

In particular, $Y$ is a boundary component of a $\Pin^+$-cobordism that is obtained by adding a 5-dimensional 2-handle to the product cobordism $M_A\times I$ 
\end{thm}

The reader is directed towards \cite{[KirbyTaylor]} and \cite[Section 2]{[Stolz]} for preliminaries and background on $\Pin^{\pm}$-structures. To distinguish the smooth structures of Theorem \ref{Theorem A} and Theorem \ref{Theorem B}, we use the $\eta$-invariant and results of Stolz \cite{[Stolz]}. The main property of this spectral invariant that we exploit is that it is a complete $\Pin^+$-cobordism invariant$\mod 2\Z$ (see Theorem \ref{Theorem Stolz}).

The choice of framing of the loops in Theorem \ref{Theorem B} is taken following Gompf \cite[Section 2]{[Gompf2]}; see Remark \ref{Remark Spheres}. The  monodromy of the mapping torus $M_A$ is described in Section \ref{Section CappellShaneson}. The first clause of Theorem \ref{Theorem B} is essentially due to Akbulut \cite[Section 4]{[Akbulut1]}; see Theorem \ref{Theorem Akbulut}. The simple construction of the 4-manifold $Y$ in terms of adding one single 2-handle resembles the construction of the Cappell-Shaneson homotopy 4-spheres; see Remark \ref{Remark Spheres}. Theorem \ref{Theorem B} is proven in Section \ref{Section Proof of Theorem B}.   

The third main result exemplifies constructions of inequivalent smooth structures obtained by twisting an embedded projective plane in the 4-manifold $Y$ of Theorem \ref{Theorem A}; see Section \ref{Section Twisting RP2} for a precise definition of this cut-and-paste operation.

\begin{thm}\label{Theorem R}There is a smoothly embedded real projective plane $\R P^2\hookrightarrow Y$ such that twisting it produces a 4-manifold $Z$ that is homeomorphic but not diffeomorphic to the circle sum $\R P^4\#_{S^1} \R P^4$ of two copies of the real projective 4-space.
There is a smoothly embedded real projective plane $\R P^2\hookrightarrow Z$ such that twisting it produces the 4-manifold $Y$ that is homeomorphic but not diffeomorphic to $S(2\gamma \oplus \R)$.
\end{thm}

To the best of our knowledge, twisting an $\R P^2$ had not been previously used to construct inequivalent smooth structures. In the orientable realm, the Dolgachev surface can be obtained from an elliptic surface $E(1) = \mathbb{CP}^2\#9\overline{\mathbb{CP}^2}$ in a similar way \cite{[AkbulutR]}; cf. \cite{[KatanagaSaekiTeragaitoYamada]}. We discuss in Section \ref{Section Proof of Theorem R} how the examples of Theorem \ref{Theorem R} are made of the same two compact pieces, and the isotopy class of the diffeomorphism used to identify their boundaries determines the homeomorphism type \cite{[CesardeSa], [KimRaymond]}. Handlebodies of the 4-manifolds of Theorem \ref{Theorem R} are drawn in Figure 4.

We enlarge the number of homeomorphism types of 4-manifolds that are known to support inequivalent smooth structures in our next main result. 

\begin{thm}\label{Theorem New Examples Structures}Let $M$ be a closed 4-manifold with $\pi_1(M) = \Z/2\times \Z/2$, Euler characteristic $\chi(M) = 1$ and such that there is an $x\in H^1(M; \mathbb{F}_2)$ with $x^4 \neq 0$. There are closed smooth homeomorphic but not diffeomorphic 4-manifolds $X_1$ and $X_2$ with the same homotopy type of $M$. 
\end{thm}

Hambleton-Hillman \cite{[HambletonHillman]} showed that the 4-manifold $M$ of Theorem \ref{Theorem New Examples Structures} is homotopy equivalent to either the non-trivial $\R P^2$-bundle over $\R P^2$ or to the quotient of the circle sum $\R P^4 \#_{S^1}\R P^4$ by a fixed-point free involution. These 4-manifolds are also obtained by twisting an embedded $\R P^2$ (see Remark \ref{Remark New Twisting}).

Next, we consider smooth 2-knots in non-orientable 4-manifolds. While embeddings of 2-spheres in orientable 4-manifolds have been widely studied in recent years, not much is known in the non-orientable case. Our next result points out further intricacies of smooth embeddings of 2-spheres inside non-orientable 4-manifolds.

\begin{thm}\label{Theorem C}There are five $\{S_i: i = 0, \ldots, 4\}$ smoothly embedded 2-spheres in $(\R P^2\times S^2)\#\mathbb{CP}^2$, which satisfy the following properties.\begin{enumerate}
\item The tubular neighborhood $\nu(S_i)$ is diffeomorphic to the $D^2$-bundle over $S^2$ with Euler number one, whose boundary is $\partial(D^2\widetilde{\times} S^2) = S^3$, for $i=0, \dots, 4$.

\item \begin{equation} S_0=\mathbb{CP}^1 \subset (\R P^2\times S^2)\#\mathbb{CP}^2 \end{equation}
and\begin{equation}
 [S_i] = [S_0]+[\{pt\} \times S^2] \in H_2((\R P^2\times S^2)\#\mathbb{CP}^2; \Z)\end{equation} for $i = 1, \ldots, 4$.
\item For $i=0, \dots, 4$, the complement of $S_i$ has cyclic fundamental group \begin{equation}\pi_1((\R P^2\times S^2)\# \mathbb{CP}^2\setminus \nu(S_i)) = \Z/2 .\end{equation}
\item The exteriors of $\{S_0, S_1, S_2\}$ are pairwise non-homotopy equivalent.
\item There are homeomorphisms of pairs\begin{equation}((\R P^2\times S^2)\# \mathbb{CP}^2, S_1)\rightarrow ((\R P^2\times S^2)\# \mathbb{CP}^2, S_3)\end{equation}and\begin{equation}((\R P^2\times S^2)\# \mathbb{CP}^2, S_2)\rightarrow ((\R P^2\times S^2)\# \mathbb{CP}^2, S_4),\end{equation}but no such diffeomorphisms exist.
\end{enumerate}
\end{thm}


Theorem \ref{Theorem C} puts on display phenomena that do not occur in the simply connected realm. Hambleton-Kreck \cite{[HambletonKreck]} and Lee-Wilczy\'nski \cite{[LeeWilczynski]} showed that for simply connected 4-manifolds any pair of locally flat 2-spheres that belong to the same homology class and whose exteriors have cyclic fundamental group are topologically isotopic; cf. \cite[Section 5]{[KasprowskiPowellRay]}. Example \ref{Example Akbulut} describes a pair of smoothly inequivalent 2- and 3-spheres in $\R P^4\#\mathbb{CP}^2$ with homeomorphic complements whose existence is an immediate consequence of work of Cappell-Shaneson \cite{[CappellShaneson]} and Akbulut \cite{[Akbulut1]}.

The paper is organized as follows. Section \ref{Section Gluck Twists} recalls the definition of Gluck twists, which are used in Section \ref{Section Bundles} to construct 2-sphere bundles over the real projective plane. In Section \ref{Section Circle sums} we explain the circle sum operation, while \ref{Section Twisting RP2} contains a description of the operation of twisting an $\R P^2$ and it includes a procedure to draw the handle diagrams of 4-manifolds that are obtained in this way. The raw materials employed in the construction of the smooth structures and the invariant used to distinguish them are described in Section \ref{Section CappellShaneson} and Section \ref{Section Invariant}, and a quick survey of some construction is given in Section \ref{Section CSRP4}. Key ingredients in the proofs of our theorems are collected in Section \ref{Section Sphere Bundle} and Section \ref{Section Result Akbulut}. The former contains a description of the 4-manifold $Y$, while in the latter we recall a theorem of Akbulut. Stabilizations with a copy of $\mathbb{CP}^2$ are discussed in Section \ref{Section CP2}.  The proofs of our main results are found in Section \ref{Section Proofs Main Results}.

\subsection{Acknowledgements} We are thankful to the referee for their thorough review, which helped us improved the manuscript. V. B. thanks Lisa Piccirillo for the nice conversation about this project. R. T. thanks Bob Gompf for sharing his knowledge through multiple e-mail exchanges. Selman Akbulut's work inspired this note, and we thank him for it. We thank Daniele Zuddas for his useful suggestions. The illustrations were done by Mateja Luke\^{z}i\^{c}. V. B. is partially supported by GNSAGA – Istituto Nazionale di Alta Matematica ‘Francesco Severi’, Italy.

\section{Background results, tools and constructions}\label{Section Tools}

\subsection{Handlebody diagrams of non-orientable 4-manifolds}\label{Section Conventions Handlebodies} Foundational results concerning handlebodies of non-orientable 4-manifolds have been obtained by Akbulut \cite{[Akbulut1], [Akbulut3]}, C\'esar de S\'a \cite{[CesardeSa]} and Miller-Naylor \cite{[MillerNaylor]}. We refer the reader to these sources for background material, and proceed to briefly recall the conventions of Akbulut \cite{[Akbulut1]}, \cite[Section 1.5]{[Akbulut3]} that we follow to draw handlebodies in this paper. The pair of round balls representing the attaching region of the orientation-reversing 1-handles are drawn with small arcs at their centers to denote that the standard oriented 1-handle identification (\cite[Figure 4.15]{[GompfStipsicz]}) is composed with a reflection across the plane orthogonal to the arcs \cite[Figure 1.25]{[Akbulut3]}. In the handlebody depicted in Figure 2, one of the round balls is at $\infty$. In this case, the identification is the composition of such reflection with the radial map that identifies a point on the sphere centered at the origin with a point on the sphere centered at $\infty$ \cite[Figure 1.25]{[Akbulut3]}.  The dotted circle notation will be used for any orientable 1-handle \cite{[Akbulut0], [CesardeSa]}.

The canonical framing of a 2-handle is the blackboard framing \cite[Section 4.5]{[GompfStipsicz]}, which depends on the chosen projection. This is the framing coming from the normal vector field to the surface of the page of the paper, and the framing of any 2-handle is taken with respect to this framing \cite[p. 117]{[GompfStipsicz]}. The notation used to record the framing of the attaching circle of a 2-handle in a non-orientable handlebody depends on whether the 2-handle runs or does not run over a non-orientable 1-handle as follows \cite[Section 0]{[Akbulut1]}.  When the 2-handle does not run over any non-orientable 1-handle, the framing is specified with a single integer number. For example, this is the case of the 0-framed circle in any of the handlebodies of Figure 1. If the 2-handle runs over a non-orientable 1-handle, an integer number is assigned to each arc of the attaching circle minus the non-orientable 1-handles. Following \cite[Section 1.5]{[Akbulut3]}, we record the framing of such a 2-handle as a pair of circled integer numbers $(m, n)\in \Z\times \Z$. For example, in the handlebody on top of Figure 1, the framing of the 2-handle that runs over the non-orientable 1-handle is $(k, 0)$; it indicates that one adds $k$ twists to the blackboard framing. Moreover, the sign of the framing of a 2-handle flips when going across an orientation-reversing 1-handle \cite[p. 78 - 79]{[Akbulut1]}.  The 3- and 4-handles do not need to be drawn to determine the diffeomorphism type of a non-orientable 4-manifold due to a foundational result of C\'esar de S\'a \cite{[CesardeSa]} and Miller-Naylor \cite{[MillerNaylor]}.

\subsection{Gluck twists}\label{Section Gluck Twists}Let $X$ be a smooth $4$-manifold and $S\hookrightarrow X$ be a smoothly embedded 2-sphere with trivial normal bundle $\nu(S) = D^2\times S^2$; there is a unique framing for $S$ \cite[Section 4.1]{[GompfStipsicz]}. The result of performing a Gluck twist to $X$ on $S$ is the 4-manifold\begin{equation}\label{Gluck Twist Manifold}X_S:= (X \backslash \nu(S)) \cup_{\varphi} (D^2\times S^2)\end{equation}obtained by using a diffeomorphism $\varphi: \partial (X \backslash \nu(S)) \rightarrow \partial(D^2\times S^2)=S^1 \times S^2$ to identify the common boundaries given by
\begin{equation}\label{Diffeo Gluck}\varphi(\theta, x) = (\theta, \rho_{\theta}(x))\end{equation} for $(\theta, x)\in S^1\times S^2$. The diffeomorphism $\rho_{\theta}: S^2\rightarrow S^2$ is a rotation of the 2-sphere through an angle $\theta\in \partial D^2 = S^1$ \cite[Section 6]{[Gluck]}. 

\begin{figure}{\label{Figure 1}}
\begin{center}
\begin{center}
\includegraphics[width=90mm,scale=0.5]{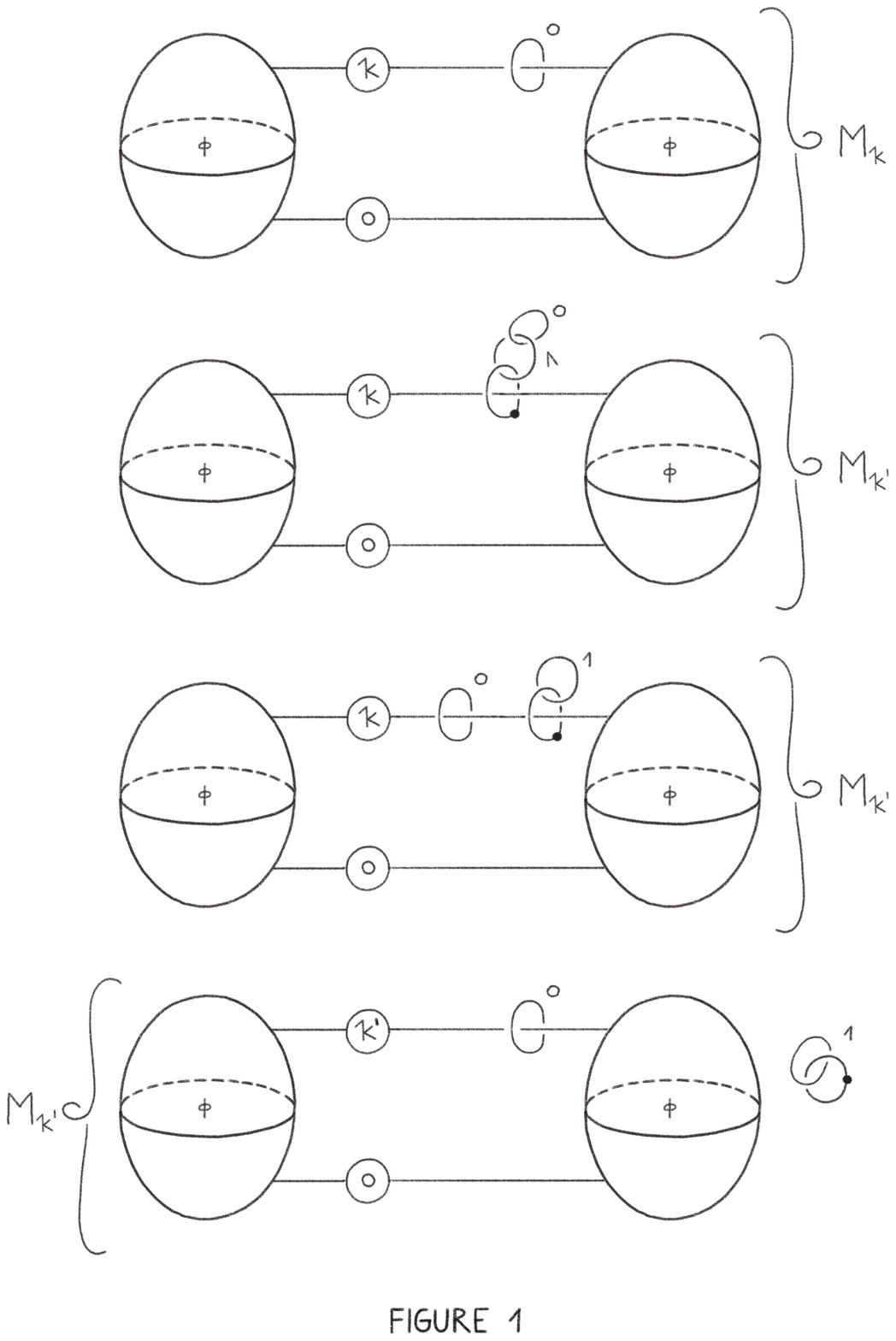}
\end{center}
\end{center}
\end{figure}

\subsection{2-sphere bundles over the real projective plane}\label{Section Bundles}There are two non-orientable total spaces of $S^2$-bundles over $\R P^2$ \cite[Chapter 12]{[Hillman]}. The non-trivial bundle is given as the quotient\begin{equation}\label{Notation Bundle}S(2\gamma \oplus \R) = \frac{S^2\times S^2}{(x, y)\sim (-x, \rot_{\pi}(y))},\end{equation} where $\rot_{\pi}: S^2\rightarrow S^2$ is a rotation of $\pi$ radians \cite{[Hillman]}. The bundle (\ref{Notation Bundle}) can also be regarded as the unit sphere bundle of the rank $3$ vector bundle over $\R P^2$ given by the Whitney sum of two copies of the canonical bundle $\gamma$ with one copy of the trivial line bundle. The latter description explains the notation (\ref{Notation Bundle}).

Both of these bundles are obtained by adding a copy of $D^2\times S^2$ to\begin{equation}\label{Basic Block Bundles}(\R P^2\times S^2)\setminus \nu (\Sigma') = Mb\times S^2 = S(2\gamma\oplus \R) \setminus \nu (\Sigma),\end{equation}where $Mb$ is the M\"obius band, $\Sigma':= \{pt\}\times S^2\subset \R P^2\times S^2$ and $\Sigma \subset S(2\gamma\oplus \R)$ is a fiber. The diffeomorphism class of the resulting 4-manifold depends on the isotopy class of the diffeomorphism $\phi: \partial Mb\times S^2\rightarrow \partial D^2\times S^2 = S^1\times S^2$ that is used to glue in $D^2\times S^2$ to (\ref{Basic Block Bundles}). Using $\phi = \id$ yields $\R P^2\times S^2$, while using the diffeomorphism (\ref{Diffeo Gluck}) yields $S(2\gamma \oplus \R)$. The next lemma is a summary of this discussion.

\begin{lemma}\label{Lemma Gluck Twist}\

$\bullet$ The non-trivial bundle $S(2\gamma\oplus \R)$ is obtained by performing a Gluck twist on $\Sigma'\subset \R P^2\times S^2$.

$\bullet$ The trivial bundle $\R P^2\times S^2$ is obtained by performing a Gluck twist on $\Sigma\subset S(2\gamma \oplus \R)$.
\end{lemma}


\begin{proof} Let $M_k$ be the 4-manifold given by the handlebody on the top level of Figure 1 along with a 3- and a 4-handle. We have $M_k = \R P^2\times S^2$ if $k = 0 \mod 2$, and $M_k = S(2\gamma\oplus \R)$ if $k = 1 \mod 2$. This can be seen as follows. There are two non-orientable 2-sphere bundles over the real projective plane \cite[Section 12.3]{[Hillman]} and both of them arise as the double of a non-orientable 2-disk bundle over $\R P^2$ \cite[Example 4.6.3]{[GompfStipsicz]}. Handlebodies of these 2-disk bundles are obtained by removing the 0-framed 2-handle, and the 3- and 4-handles from the handlebody of $M_K$ \cite[p. 79]{[Akbulut1]} (cf. \cite{[MillerNaylor]}). The 0-framed 2-handle (coming from the double construction) in the handlebody of $M_K$ of Figure 1 represents the tubular neighborhood $\nu(S) = D^2\times S^2$ of an embedded 2-sphere $S\subset M_k$ \cite[Exercise 5.4.3. (d)]{[GompfStipsicz]}. Let $M_{k'}$ be the 4-manifold obtained by doing a Gluck twist on this 2-sphere. Its handlebody is drawn on the second stage from top to bottom in Figure 1 following Akbulut-Yasui's method \cite[Figure 1]{[AkbulutYasui]}. Take the 0-framed 2-handle that links the 1-framed 2-handle, and slide it under the 1-handle represented by the dotted circle to obtain the third stage of Figure 1. Sliding the $k$-framed 2-handle over the 1-framed 2-handle and then under the dotted circle results in the handlebody drawn at the bottom of Figure 1: the framing $k'$ is an odd integer if and only if $k$ is an even integer. The bottom level of Figure 1 is the handlebody of $M_k'$, where the Hopf link of the dotted circle and the 1-framed unknot is a cancelling pair of 1- and 2-handles. 


\end{proof}

\subsection{Circle sums}\label{Section Circle sums}
Let $\{X_i: i=1,2 \}$ be closed 4-manifolds with $H^1(X_i; \mathbb{Z}/2) \cong \mathbb{Z}/2$ and that admit a $\text{Pin}^{+}$-structure, and let $D^3 \widetilde{\times} S^1$ be the non-trivial $3$-disk bundle over $S^1$ endowed with a fixed $\text{Pin}^{+}$-structure. There is a diffeomorphism $\nu(\gamma_i)\rightarrow D^3\widetilde{\times} S^1$, where $\nu(\gamma_i)$ is the tubular neighborhood of an orientation-reversing simple loop $\gamma_i\subset X_i$. The circle sum of $X_1$ and $X_2$ along $\gamma_1$ and $\gamma_2$ (cf. \cite{[HambletonKreckTeichner]}) is defined as the smooth $\text{Pin}^{+}$ $4$-manifold\begin{equation}\label{circle}X_1\#_{S^1} X_2 :=(X_1 \setminus \nu(i_1(D^3 \widetilde{\times} S^1)) \cup_{\varphi} (X_2 \setminus \nu(i_2(D^3 \widetilde{\times} S^1))\end{equation}where the gluing diffeomorphism is\begin{equation}\label{Gluing CS}\varphi:= i_2|_{\partial D^3 \widetilde{\times} S^1} \cdot (i_1|_{\partial D^3 \widetilde{\times} S^1})^{-1}: \partial(X_1 \setminus \nu(i_1(D^3\widetilde{\times} S^1)) \to \partial (X_2 \setminus \nu(i_2(D^3\widetilde{\times} S^1))\end{equation}and $i_i: D^3\widetilde{\times} S^1 \to X_i$ are embeddings representing the loop generating $H^1(X_i; \mathbb{Z}/2)$ via the universal coefficient theorem for $i=1,2$. Moreover, we ask $i_1$ to preserve the $\text{Pin}^{+}$-structure and $i_2$ to reverse it. 

\begin{lemma}\label{cobordism}Let $X_1$ and $X_2$ be two closed smooth non-orientable 4-manifolds with $H^1(X_i; \Z/2) = \Z/2$ that admit a $\Pin^+$-structure and  let $\gamma_i\subset X_i$ be an orientation-reversing simple loop for $i = 1, 2$. There is a $\text{Pin}^+$-cobordism between the disjoint union $X_1 \sqcup X_2$ and their circle sum (\ref{circle}). 
\end{lemma}

\begin{proof}(cf. \cite[p. 160]{[Stolz]}). A cobordism $W$ is obtained from $(X_1\sqcup X_2)\times [0, 1]$ by identifying its upper lid with the boundary of $(S^2\widetilde{\times} S^1)\times [0, 1]$ as\begin{equation}\label{Pin Cobordism Lemma}(X_1\setminus \nu(i_1(D^3\widetilde{\times} S^1))\times \{1\})\cup ((S^2\widetilde{\times} S^1)\times [0, 1])\cup (X_2\setminus\nu(i_2(D^3\widetilde{\times} S^1))\times \{1\})\end{equation}using the map (\ref{Gluing CS}). A $\Pin^+$-structure on $W$ is obtained by gluing together the matching $\Pin^+$-structures on the building blocks in (\ref{Pin Cobordism Lemma}) and it induces a $\Pin^+$-structure on $X_1\#_{S^1}X_2$ \cite{[Kirby], [KirbyTaylor]}.
\end{proof}

We point out that the circle sum of two irreducible 4-manifolds need not result in an irreducible 4-manifold.

\begin{lemma}\label{Lemma CS Reducible} The circle sum $S(2\gamma \oplus \R)\#_{S^1}\R P^4$ is diffeomorphic to $\R P^4\#(S^2\times S^2)$.  
\end{lemma}

\begin{proof}A handlebody of $\R P^4\#(S^2\times S^2)$ is obtained by modifying the handlebody on the third stage from top to bottom of Figure 1 as follows: erase the 0-framed circle, substitute the dotted circle with a 2-framed circle, change the framing of the linking circle from 1 to 0, and set $k = 1$. The 3-handle and the 4-handle do not need to be drawn \cite{[CesardeSa], [MillerNaylor]}. A handlebody for $S(2\gamma\oplus \R)\#_{S^1} \R P^4$ is obtained by adding a $(1, 0)$-framed 2-handle to the first stage from bottom to top of Figure 3. Slide the $(-1, 0)$-framed 2-handle over the $(1, 0)$-handle and then slide the latter over the other $(1, 0)$-handle to obtain the handlebody of $\R P^4\#(S^2\times S^2)$.


\end{proof}

\subsection{Twisting an embedded $\R P^2$}\label{Section Twisting RP2}Let $X$ be a smooth $4$-manifold and let $R\hookrightarrow X$ be a smoothly embedded real projective plane with tubular neighborhood $\nu(R)$ diffeomorphic to the non-trivial and non-orientable 2-disk bundle\begin{equation}D^2\widetilde{\times} \R P^2 = \frac{D^2\times S^2}{(d, y)\sim (-d, - y)}\end{equation} for every $(d, y)\in D^2\times S^2$ \cite[Section 5]{[HambletonHillman]}. The result of twisting $\nu(R)\subset X$ is the 4-manifold\begin{equation}\label{Manifold Surgery RP2}X_R: = (X\setminus \nu(R))\cup_{\varphi'}(D^2\widetilde{\times} \R P^2)\end{equation}obtained by using a diffeomorphism $\varphi': \partial (X \backslash \nu(R)) \rightarrow \partial(D^2\widetilde{\times} \R P^2) = S^2\widetilde{\times} S^1$ to identify the common boundaries given by\begin{equation}\label{Diffeo NGluck}\varphi'(x, \theta) = (\rho_{\theta}(x), \theta)\end{equation} for $\theta$ in the circle base and $x$ in the 2-sphere fiber of the 2-sphere bundle over the circle $S^2\widetilde{\times} S^1$ \cite[Section 2]{[CesardeSa]}, \cite{[KimRaymond]}.

 \begin{figure}{\label{Figure 1}}
\begin{center}
\includegraphics[width=100mm,scale=0.2]{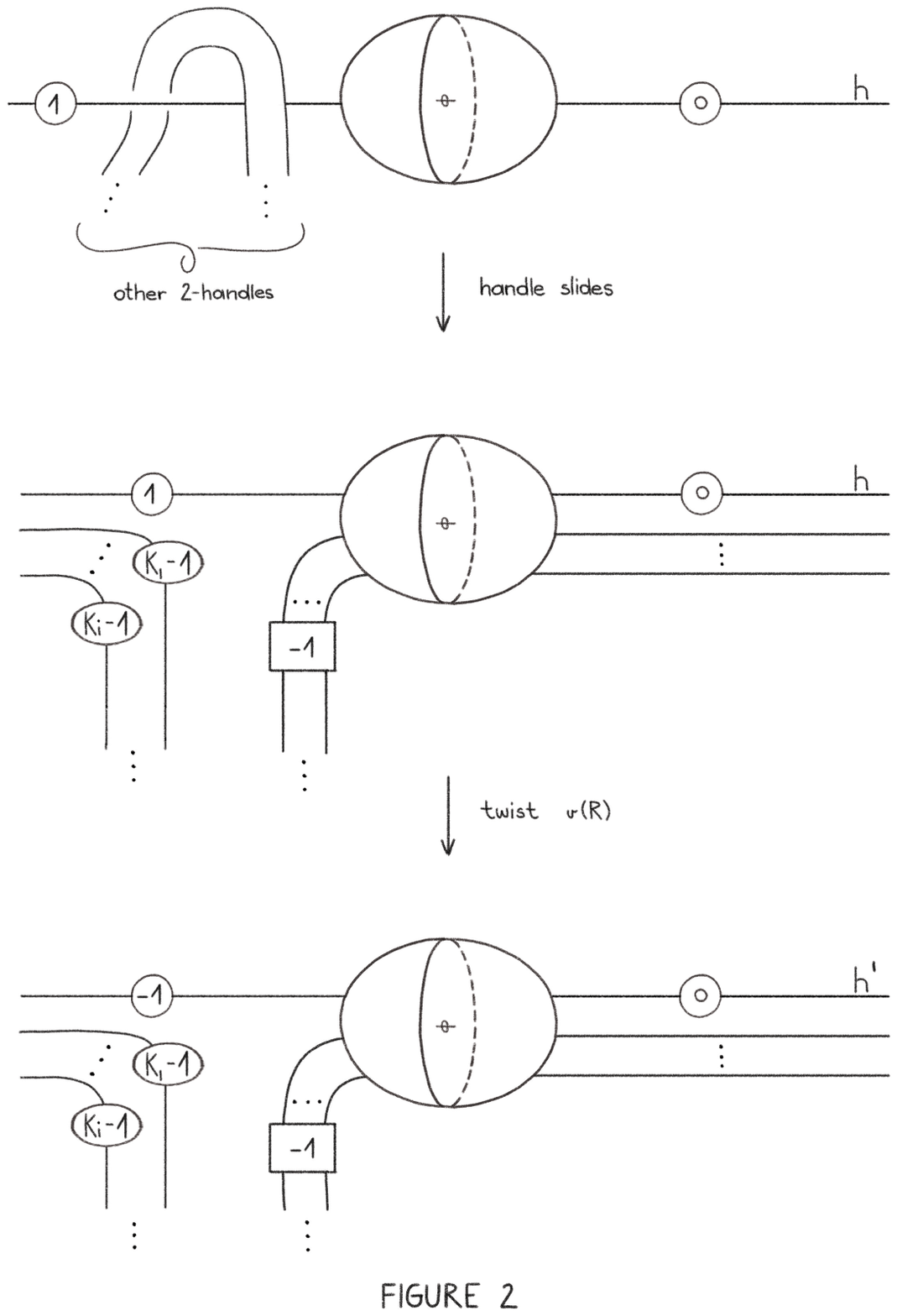}
\end{center}
\end{figure}

We describe a method to draw a handlebody of the construction (\ref{Manifold Surgery RP2}) by building heavily on the operation of blowing down an $\R P^2$ that was introduced by Akbulut in \cite[p. 79 - 80]{[Akbulut1]}. The method is described in Figure 2, where one of the round balls of the attaching region of the 1-handle is at the point of infinity and a single round ball is drawn. The handlebody of $X$ is depicted by the diagram at the top level of Figure 2, where the union of the 1-handle and the 2-handle $h$ with strands drawn horizontally with framing $(1, 0)$ is the tubular neighborhood $\nu(R) = D^2\widetilde{\times} \R P^2$. We first argue that this is the only case that needs to be considered. If the 2-handle in the handle decomposition of $\nu(R)$ has framing $(-1, 0)$, first rotate the round ball drawn 360 degrees around the axis indicated by the arc (\cite[p. 79]{[Akbulut1]}) to change its framing to $(0, -1)$, and then transfer the framing through the orientation-reversing 1-handle to obtain a 2-handle with framing $(1, 0)$ as in the top level of Figure 2 (see \cite[Figure 1.26]{[Akbulut3]}). Denote the 2-handles that link $h$ by $\{K_1, \ldots, K_i\}$ with their respective framings given by $\{k_1, \ldots, k_i\}$. Slide the attaching circles of all the 2-handles $\{K_r: r = 1, \ldots, i\}$ over $h$ until they do not link it anymore and obtain the second level of Figure 2 (cf. \cite[p. 80]{[Akbulut1]}). The framings are now $\{k_1 - 1, \ldots, k_i - 1\}$ and the strands have a $-1$ full twist as indicated in the handlebody of $X$ drawn in the second level of Figure 2. At this point we perform the twist (\ref{Manifold Surgery RP2}) on $\nu(R)$ and draw the handlebody of $X_R$ at the bottom level of Figure 2. This is done by first blowing down the real projective plane $R$ by erasing the 2-handle $h$ from the handlebody of $X$ drawn in the second level of Figure 2, and then adding a 2-handle $h'$ with framing $(-1, 0)$ that replaces $h$. 

So far we have described a method that produces a 4-manifold $X'$, which is obtained by carving $\nu(R)$ out of $X$ and gluing it back in using a self-diffeomorphism of the boundary $\partial(X\setminus \partial \nu(R)) = S^2\widetilde{\times} S^1$. To conclude that $X'$ is diffeomorphic to (\ref{Manifold Surgery RP2}), we argue as follows  (cf. \cite[Proof Lemma 2.2]{[AkbulutYasui]}). A result of Kim-Raymond \cite{[KimRaymond]} on the diffeotopy group of $S^2\widetilde{\times} S^1$ (cf. \cite[Remark 3.12]{[MillerNaylor]}) implies that we need to check that the gluing diffeomorphism (\ref{Diffeo NGluck}) does not extend to a self-diffeomorphism of $D^2\widetilde{\times}\R P^2$. This is done in the following example. 

\begin{figure}{\label{Figure 1}}
\begin{center}
\includegraphics[width=100mm,scale=0.4]{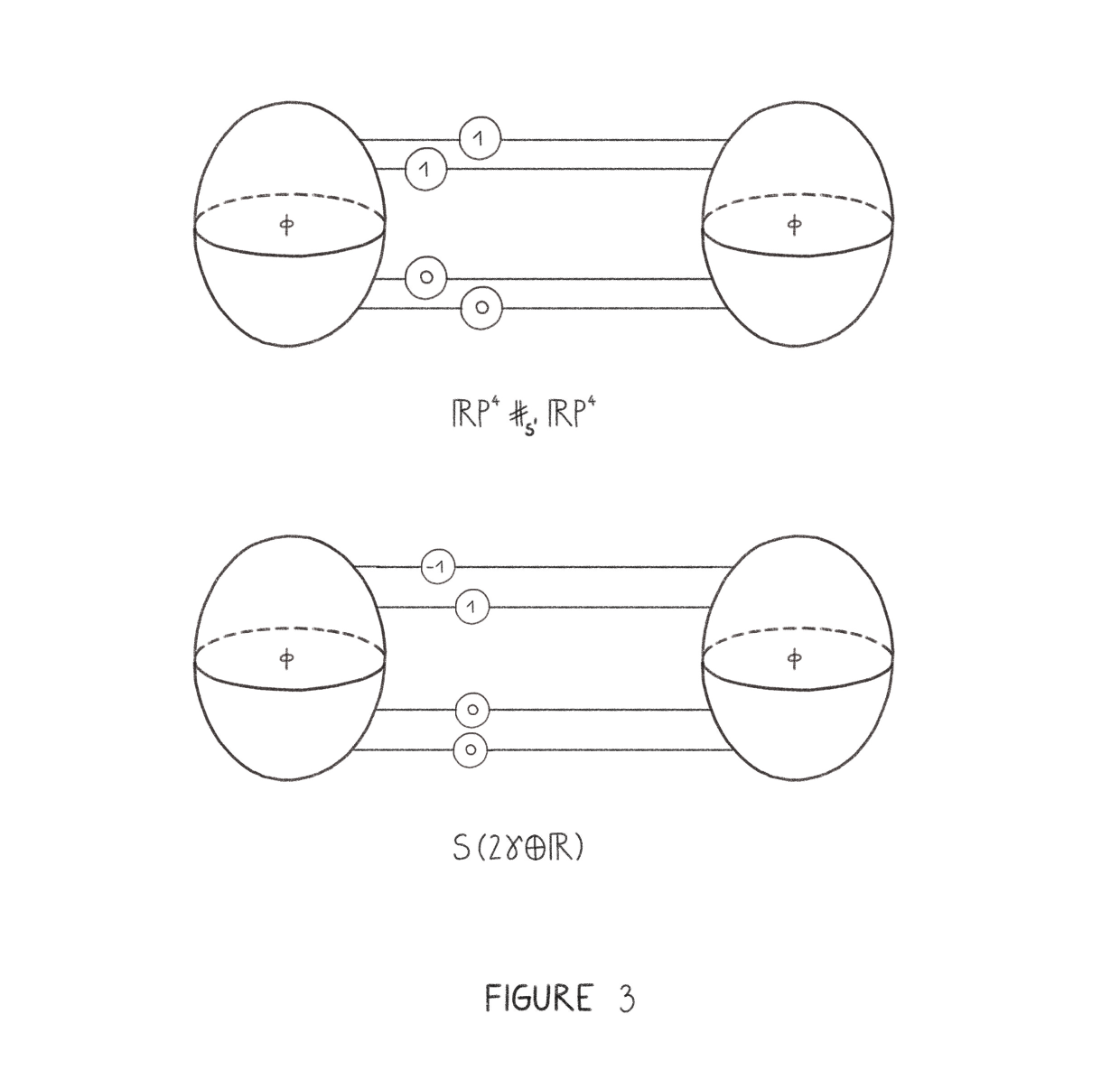}
\end{center}
\end{figure}


\begin{example}\label{Example Basic}The homotopy equivalence classes of the double\begin{equation}\label{S_3}S(2\gamma\oplus \R) = (D^2\widetilde{\times} \R P^2)\cup_{\id}(D^2\widetilde{\times} \R P^2)\end{equation}and of the circle sum\begin{equation}\label{S_4}\R P^4\#_{S^1} \R P^4 = (D^2\widetilde{\times} \R P^2)\cup_{\varphi'}(D^2\widetilde{\times} \R P^2)\end{equation}do not coincide \cite{[KimKojimaRaymond]}. Handlebodies of these two non-orientable 4-manifolds are drawn in Figure 3 with an additional 3- and 4-handle in each, making visible the procedure of twisting an $\R P^2$.
\end{example}

Hambleton-Hillman \cite{[HambletonHillman]} studied the homotopy and homeomorphism types of closed 4-manifolds whose universal cover is $S^2\times S^2$. A representative of a homotopy equivalence class in the quadratic 2-type of $PD_4$-complexes with Euler characteristic one and fundamental group $\Z/2\times \Z/2$ is a 4-manifold $W$, which arises as the quotient space of $\R P^4\#_{S^1} \R P^4$ by a free involution \cite[Definition 7.1]{[HambletonHillman]}. In order to exemplify the utility of the construction (\ref{Manifold Surgery RP2}), we now remark that $W$ can be obtained from the non-trivial bundle $\R P^2\widetilde{\times} \R P^2$ by twisting an $\R P^2$.

\begin{example}\label{Example New Construction}Use the decomposition $\R P^2 = D^2\cup Mb$ to deconstruct the bundle as\begin{equation}\label{Decomposition Bundle}\R P^2\widetilde{\times} \R P^2 = (D^2\widetilde{\times} \R P^2) \cup_{\varphi} (Mb\widetilde{\times} \R P^2),\end{equation} where $Mb\widetilde{\times} \R P^2$ is a non-trivial M\"obius band bundle over $\R P^2$ with boundary $\partial (Mb\widetilde{\times} \R P^2) = S^2\widetilde{\times} S^1$. In the previous equality, there is a choice of diffeomorphism $f:\partial(Mb\widetilde{\times} \R P^2)\rightarrow S^2\widetilde{\times} S^1$ that is made as follows. A $\Pin^+$-structure on $\R P^2\widetilde{\times} \R P^2$ induces a $\Pin^+$-structure on each of the building blocks in (\ref{Decomposition Bundle}) as codimension zero submanifolds \cite[Construction 1.5]{[KirbyTaylor]}. The induced $\Pin^+$-structures restrict to the same $\Pin^+$-structure $\phi^{S^2\widetilde{\times} S^1}$ on each of the $S^2\widetilde{\times} S^1$ boundary components $\partial(D^2\widetilde{\times} \R P^2)$ and $\partial(Mb\widetilde{\times} \R P^2)$. The diffeomorphism $f$ is chosen so that it preserves the bundle structure and satisfies $f^\ast \phi^{S^2\widetilde{\times} S^1} = \phi^{S^2\widetilde{\times} S^1}$. A result of Kim-Raymond \cite{[KimRaymond]} says that there are two isotopy classes of self-diffeomorphisms of $S^2\widetilde{\times} S^1$ that respect the bundle structure. Taking into account the induced $\Pin^+$-structures, we conclude that the gluing diffeomorphism $\varphi$ in (\ref{Decomposition Bundle}) is isotopic to the identity map. From (\ref{Decomposition Bundle}), we see that $S(2\gamma \oplus \R)$ is obtained by taking a double cover (fiberwise) of $\R P^2\widetilde{\times} \R P^2$. Using a self-diffeomorphism $\varphi'$ of $S^2\widetilde{\times} S^1$ that preserves the bundle structure and that is not isotopic to the identity map, we obtain\begin{equation}\label{Decomposition Not Bundle}W = (D^2\widetilde{\times} \R P^2) \cup_{\varphi'} (Mb\widetilde{\times} \R P^2)\end{equation} with $\R P^4\#_{S^1}\R P^4$ as a double cover. 
\end{example}

The 4-manifolds $X$ and $X_R$ are $\mathbb{CP}^2$-stable in the following sense.

\begin{lemma}\label{Lemma Stabilize Twisting RP2}There is a diffeomorphism\begin{equation}\label{Diffeomorphism Lemma}X\# \mathbb{CP}^2 \rightarrow X_R\# \mathbb{CP}^2.\end{equation} 
\end{lemma}

\begin{proof}Erase the 2-handle $h$ from the handlebody in the middle level of Figure 2 to produce the compact 4-manifold $X_0' = X\setminus \nu(R)$. To draw a handlebody of $X\#\overline{\mathbb{CP}^2}$, add a -1-framed unknot to the middle diagram of Figure 2, and slide the 2-handle $h$ over it. The new 2-handle $h$ now has $(0, 0)$-framing and it is linked by a -1-framed circle. This exhibits a decomposition of the connected sum as\begin{equation}\label{Decomposition L}X\#\overline{\mathbb{CP}^2} = X_0'\cup ((D^2\widetilde{\times}\R P^2)\#\overline{\mathbb{CP}^2}).\end{equation}To see that there is a decomposition\begin{equation}\label{Decomposition R}X_R\#\mathbb{CP}^2 = X_0'\cup ((D^2\widetilde{\times}\R P^2)\#\mathbb{CP}^2),\end{equation}draw a handlebody for this connected sum by adding a 1-framed unknot to the bottom diagram of Figure 2 and sliding the $(-1, 0)$-framed 2-handle $h'$ over it. These two handlebodies of $X\#\overline{\mathbb{CP}^2}$ and $X_R\#\mathbb{CP}^2$ differ solely on the sign of the circle that links the $(0, 0)$-framed 2-handle and, in particular, the gluing maps in (\ref{Decomposition L}) and (\ref{Decomposition R}) coincide. Since $M\# \mathbb{CP}^2$ is diffeomorphic to $M\#\overline{\mathbb{CP}^2}$ for every non-orientable 4-manifold $M$, we conclude that the diffeomorphism (\ref{Diffeomorphism Lemma}) exists.

\end{proof}

There are instances where twisting an $\R P^2$ does not change the diffeomorphism type of $X$ and it is equivalent to twisting a non-orientable 1-handle cf. \cite[p. 82]{[Akbulut1]}, \cite[Figure 5.42]{[GompfStipsicz]}. In such a  scenario, $X$ is diffeomorphic to $X_R$ and the result of this self-diffeomorphism is to relate the $\Pin^+$-structures on $X$ as we now exemplify. 

\begin{example}\label{Example RP4}Deconstruct the real projective $4$-space as\begin{equation}\label{RP4}\R P^4 = (D^2\widetilde{\times}\R P^2)\cup (D^3\widetilde{\times} S^1)\end{equation}and twist the obvious $\R P^2$. Since any self-diffeomorphism of $S^2\widetilde{\times} S^1$ extends to a self-diffeomorphism of $D^3\widetilde{\times} S^1$, twisting the $\R P^2$ in (\ref{RP4}) does not change the diffeomorphism type. There are two $\Pin^+$-structures $(\R P^4, \pm \phi_{\R P^4})$ and they are mutual inverses in the fourth $\Pin^+$-cobordism group $\Omega_4^{\Pin^+} = \Z/16$ \cite{[KirbyTaylor]}. The twist of the non-orientable 1-handle gives a self-diffeomorphism of $\R P^4$ that maps $\pm \phi_{\R P^4}$ to $\mp \phi_{\R P^4}$. We use the following convention to distinguish the $\Pin^+$-structures:\begin{equation}\eta(\R P^4, \pm \phi_{\R P^4}) = \pm \frac{1}{8} \mod 2\Z.\end{equation}
\end{example}

Although this cut-and-paste operation is of a local nature in the sense that $X\setminus D^2\widetilde{\times} \R P^2\subset X_R$ is a codimension zero submanifold, the topology of $X_R$ can well be drastically different to that of $X$ as exemplified in the following lemma. 

\begin{lemma}Twisting an $\R P^2$ need not preserve irreducibility.
\end{lemma}

\begin{proof} The circle sum of three copies of the real projective 4-space $\#_{S^1}3\cdot \R P^4$ contains a copy of $D^2\widetilde{\times} \R P^2$. Twisting this $\R P^2$ yields $S(2\gamma \oplus \R)\#_{S^1}\R P^4$ (cf. Figure 3), which is diffeomorphic to $\R P^4\# (S^2\times S^2)$ by Lemma \ref{Lemma CS Reducible}. 
\end{proof}

\subsection{$\Pin^{\pm}$-structures on a mapping torus and framings on its cross-section}\label{Section CappellShaneson}As it was done by Cappell-Shaneson in \cite[Section 2]{[CappellShaneson]}, consider the diffeomorphism\begin{equation}\label{Choice Diffeo}\varphi_{A}: T^3\rightarrow T^3 = \R^3/\Z^3\end{equation} whose induced map $(\varphi_{A})_{\ast}: H_1(T^3; \Z)\rightarrow H_1(T^3; \Z) = \Z^3$ is given by the matrix\begin{equation}\label{Matrix 1}A =\left( \begin{array}{ccc}
0 & 1  & 0   \\
0 & 0  & 1  \\
- 1 & 1  & 0
\end{array} \right)\in \Gl(3; \Z).\end{equation}The mapping torus\begin{equation}\label{Mapping Torus}M_A := (T^3\times [0, 1])/\sim\end{equation} of the diffeomorphism (\ref{Choice Diffeo}) under the identification $(x, 0) \sim (\varphi_A(x), 1)$ is the total space of a non-orientable $T^3$-bundle over $S^1$. Let $\nu_0$ be the tubular neighborhood of the loop\begin{equation}\label{Loop Section}\alpha_0:= ((\{(x_0\} \times [0, 1])/\sim) \subset M_A,\end{equation} where $\varphi_A(x_0) = x_0$ for $x_0 \in T^3$. The loop $\alpha_0$ inherits a canonical framing from $\{x_0\}\times [0, 1]$ \cite{[AitchisonRubinstein]}, \cite[Section 2]{[Gompf2]}. We will call the other framing the non-canonical framing following Gompf \cite[p. 1668]{[Gompf2]}. The compact 4-manifold $\nu_0$ is diffeomorphic to $D^3\widetilde{\times} S^1$ with $\partial \nu_0 = S^2\widetilde{\times} S^1$.

The 4-manifold (\ref{Mapping Torus}) admits a pair of $\Pin^{\dagger}$-structures for $\dagger = +, -$ and we record some of their properties in the following proposition.

\begin{proposition}\label{Proposition Pin Mapping}\

$\bullet$ There is a pair of $\Pin^{\dagger}$-structures $(M_A, \pm \phi^{\dagger})$ for $\dagger = +, -$. 

$\bullet$ The $\Pin^-$-structure $(M_A, \phi^-)$ induces the $\Pin^-$-structure on $M_A\setminus \nu_0$ that restricts to the $\Pin^-$-structure on $\alpha_0$ that extends to the 2-disk $D^2$. 

$\bullet$ The $\Pin^+$-structure $(M_A, \phi^+)$ induces the $\Pin^+$-structure on $M_A\setminus \nu_0$ that restricts to the $\Pin^+$-structure on $\alpha_0$ that extends to the M\"obius band $Mb$. 

\end{proposition}

The proof of Proposition \ref{Proposition Pin Mapping} follows from well-known results \cite{[KirbyTaylor], [Stolz]}. We will write $(M_A, \phi^+) = (M_A, \phi^{M_A})$ for the $\Pin^+$-structure on $M_A$ in the sequel. 

The circle sum $M\#_{S^1}M_A$ along an orientation-reversing simple loop in a 4-manifold $M$ and the loop (\ref{Loop Section}) is a potential exotic smooth structure on $M$. The homeomorphism type of these 4-manifolds is the same by the following result. 

\begin{theorem}\label{Theorem CSF}Cappell-Shaneson \cite{[CappellShaneson]}, Freedman \cite{[Freedman]}. Let $M$ be a closed smooth non-orientable 4-manifold and let $\alpha\subset M$ be a simple loop whose homotopy class generates $\pi_1(M) = \Z/2$. Let $\alpha_0\subset M_A$ be the orientation-reversing simple loop in the mapping torus $M_A$ that was defined in (\ref{Loop Section}). The 4-manifolds $M$ and $M\#_{S^1} M_A = (M\setminus \nu(\alpha)) \cup (M\setminus \nu_0)$ are homeomorphic.
\end{theorem}

A proof of Theorem \ref{Theorem CSF} is obtained by generalizing an argument due to Ruberman \cite[Proof of Theorem 2]{[Ruberman]}. It also follows as a corollary of work of Hambleton-Kreck-Teichner \cite[Theorem 1]{[HambletonKreckTeichner]}. We now provide a sketch of the argument for the convenience of the reader (cf. \cite[p. 159]{[Stolz]}). A result of Cappell-Shaneson \cite[Theorem 3.1]{[CappellShaneson]} says that $M$ is simple homotopy equivalent to the circle sum $M\#_{S^1} M_A$ via a simple homotopy equivalence whose normal invariant is the non-trivial element in the kernel of $[M, G/O]\rightarrow [M, G/TOP]$. Using the surgery exact sequence \cite[Theorems 10.3 and 10.5]{[Wall]}\begin{equation}\rightarrow L_5(\Z[\pi_1(M)], w_1(M))\rightarrow \mathcal{S}^{TOP}(M)\rightarrow [M, G/TOP]\rightarrow\end{equation} one concludes that there is a topological s-cobordism between these two 4-manifolds. Given that $\pi_1(M) = \Z/2$, the surgery group vanishes (\cite[Section 4]{[KasprowskiPowellRay1]}) and results of Freedman imply that $M$ is homeomorphic to $M\#_{S^1}M_A$ \cite[Chapter 11]{[FreedmanQuinn]} \cite{[Freedman]}.

\subsection{An invariant to distinguish smooth structures}\label{Section Invariant} We collect results of Stolz \cite{[Stolz]} that are employed in this note to distinguish two smooth structures. Let $(M, \phi^M)$ be a closed smooth 4-manifold with $\Pin^+$-structure $\phi^M$ and let $\eta(M, \phi^M)$ be the $\eta$-invariant of the twisted Dirac operator on $M$ \cite[Section 2]{[Stolz]}.

\begin{theorem}\label{Theorem Stolz}Stolz \cite[Propositions 4.3 and 7.3]{[Stolz]}. Let $(M, \phi^M)$ be a closed smooth 4-manifold with $\Pin^+$-structure $\phi^M$.

$\bullet$ The invariant $\eta(M, \phi^M) \mod 2\Z$ is a complete $\Pin^+$-bordism invariant.

$\bullet$ The value of the $\eta$-invariant of the mapping torus (\ref{Mapping Torus}) is\begin{equation}\eta(M_A, \pm \phi^{M_A}) = 1 \mod 2\Z\end{equation} for both $\Pin^+$-structures $(M_A, \pm \phi^{M_A})$.

$\bullet$ Let $\alpha_1\subset M_1$ be an orientation-reversing simple loop in a closed 4-manifold $M_1$ with $\Pin^+$-structure $\phi^{M_1}$, which restricts to the $\Pin^+$ structure on the submanifold $\alpha_1$ that extends to the 2-disk. There is a $\Pin^+$-structure $\phi^M$ on the circle sum\begin{equation}\label{New Manifold}M:= M_1\#_{S^1}M_A = (M_1\setminus \nu(\alpha_1))\cup (M_A\setminus \nu_0)\end{equation} for which the $\eta$-invariant is\begin{equation}\eta(M, \phi^M)  = \eta(M_1, \phi^{M_1}) + 1 \mod 2\Z.\end{equation}
\end{theorem}

A consequence of Theorem \ref{Theorem Stolz} is the following result. We say that two smooth 4-manifolds $X$ and $Y$ are stably diffeomorphic if there is an $n\in \N$ such that the connected sums $X\#n(S^2\times S^2)$ and $Y\#n(S^2\times S^2)$ are diffeomorphic.

\begin{theorem}\label{Theorem Stolz Main} Stolz  \cite[Theorems B and 7.4]{[Stolz]}. Let $\{M_i: i = 1, 2\}$ be closed smooth non-orientable 4-manifolds that admit a $\Pin^+$-structure $\phi^{M_i}$ for $i = 1, 2$. Suppose $H^1(M_i; \Z/2) = \Z/2$ for $i = 1, 2$. If\begin{equation}\label{Different Eta}\eta(M_1, \phi^{M_1})\neq \pm \eta(M_2, \phi^{M_2}) \mod 2\Z,\end{equation} then $M_1$ and $M_2$ are neither diffeomorphic nor stably diffeomorphic. 
\end{theorem}


\begin{remark}\label{Remark Theorem Stolz Main}Stolz showed that the condition (\ref{Different Eta}) is an obstruction to the existence of a diffeomorphism $f:M_1\rightarrow M_2$ such that $f^{\ast}\phi^{M_2} = \phi^{M_1}$. The hypothesis $H^1(M_i; \Z/2) = \Z/2$ of Theorem \ref{Theorem Stolz Main} implies that there are only two possible choices of $\Pin^+$-structures on $M_i$ for $i = 1, 2$ \cite[Corollary 6.4]{[Stolz]} and these are mutual inverse elements in the fourth $\Pin^+$-cobordism group \cite[p. 190]{[KirbyTaylor]}. 

\end{remark}

We now record the changes of the $\eta$-invariant under the cut-and-paste constructions discussed in Section \ref{Section Twisting RP2}. Let $X_0 = (X\setminus \nu(R))\cup (D^3\widetilde{\times} S^1)$ be the 4-manifold obtained by blowing down a smoothly embedded real projective plane $R\subset X$ as in \cite[p. 94]{[Akbulut1]} and let $X_R$ be the 4-manifold obtained by twisting it as in Section \ref{Section Twisting RP2}. 

\begin{proposition}\label{Proposition Changes Eta}Let $X$ be a closed smooth 4-manifold that contains a smoothly embedded real projective plane $R\subset X$ with tubular neighborhood $\nu(R) = D^2\widetilde{\times} \R P^2$. Suppose there is a $\Pin^+$-structure $(X, \phi^X)$ that restricts to a $\Pin^+$-structure on $(\nu(R), \phi^{\nu(R)})$ that extends to the $\Pin^+$-structure $(\R P^4, \pm \phi^{\R P^4})$ with the convention introduced in Example \ref{Example RP4}. There are $\Pin^+$-structures $(X_R, \phi^{X_R})$ and $(X_0, \phi^{X_0})$ such that the identities\begin{equation}\label{Value Eta B1}\eta(X_0, \phi^{X_0}) = \eta(X, \phi^X) \mp \frac{1}{8} \mod 2\Z,\end{equation}\begin{equation}\label{Value Eta B1'}\eta(X_0, \phi^{X_0}) = \eta(X_R, \phi^{X_R}) \pm \frac{1}{8} \mod 2\Z\end{equation}and\begin{equation}\label{Value Eta B2}\eta(X_R, \phi^{X_R}) = \eta(X, \phi^X) \mp \frac{2}{8} \mod 2\Z\end{equation}hold.
\end{proposition}

\begin{proof} Assemble together $(X_0, \phi^{X_0})$ and $(X_R, \phi^{X_R})$ from $\Pin^+$-structures on the building blocks that match on the common boundaries. The $\Pin^+$-structure $\phi^X$ induces $\Pin^+$-structures $(X\setminus \nu(R), \phi^X|)$ and $(\nu(R), \phi^{\nu(R)})$ since these are codimension zero submanifolds of $X$ \cite{[KirbyTaylor]}, and both these $\Pin^+$-structures coincide on $\partial(X\setminus \nu(R)) = \partial \nu(R)$. Recall that the self-diffeomorphism (\ref{Diffeo NGluck}) of $S^2\widetilde{\times} S^1$ identifies the two possible choices of  $\Pin^+$-structures cf. \cite[p. 35]{[Kirby]}. The existence of the $\Pin^+$-structure $(X_R, \phi^{X_R})$ follows. Moreover, either choice of $\Pin^+$-structure $(D^3\widetilde{\times} S^1, \pm \phi)$ yields a $\Pin^+$-structure $(X_0, \phi^{X_0})$. To see that the identities (\ref{Value Eta B1}) and (\ref{Value Eta B1'}) hold, one constructs a $\Pin^+$-cobordism as in the proof of Lemma \ref{cobordism}. Such cobordism is obtained by patching together $X_0\times [0, 1]$ and $\R P^4\times [0, 1]$ along $X_0\times \{1\}$ and $\R P^4\times \{1\}$. We obtain $\Pin^+$ cobordisms $(V; X, X_0\sqcup \R P^4)$ and $(V_R; X_R, X_0\sqcup \R P^4)$ inducing  different $\Pin^+$-structures on $\R P^4$. The $\eta$-invariant is a complete $\Pin^+$-cobordism invariant and it is additive with respect to disjoint unions. It follows that (\ref{Value Eta B1}) and (\ref{Value Eta B1'}) hold. These two values imply (\ref{Value Eta B2}).

\end{proof}

\subsection{An inequivalent smooth structure and an embedded 2-sphere}\label{Section Sphere Bundle} With the purpose of unveiling some of its properties, we now describe in detail the construction of the 4-manifold $Y$ of Theorem \ref{Theorem A} given in \cite{[Torres]}. 

\begin{proposition}\label{Proposition Sphere Bundle}(cf. \cite[Section 3.2]{[Torres]}). There is a smooth 4-manifold $Y$ that satisfies the following properties.\begin{itemize}
\item It is homeomorphic to $S(2\gamma\oplus \R)$.
\item The value of its $\eta$-invariant is\begin{equation}\label{Eta Sphere}\eta(Y, \phi^Y_i) = 1 \mod 2\Z\end{equation} for both $\Pin^+$-structures $\{(Y, \phi^Y_i): i = 1, 2\}$, and $Y$ is neither diffeomorphic nor stably diffeomorphic to $S(2\gamma \oplus \R)$. 
\item There is a smoothly embedded 2-sphere $\Sigma\subset Y$ with trivial tubular neighborhood $\nu (\Sigma) = D^2\times S^2$ and non-trivial homology class $[\Sigma] \neq 0 \in H_2(Y; \Z/2)$.

\end{itemize}
\end{proposition}

The computation of (\ref{Eta Sphere}) appeared in \cite[Proposition 4, Section 3.2]{[Torres]}. We provide the argument in order to make this paper more self-contained.

\begin{proof}The 4-manifold in the statement of the proposition is the circle sum\begin{equation}\label{Construction Y}Y = S(2\gamma\oplus \R)\#_{S^1} M_A = (S(2\gamma\oplus \R)\backslash \nu(\alpha))\cup (M_A\backslash \nu_0),\end{equation}where $\alpha\subset S(2\gamma\oplus \R)$ is a loop whose homotopy class generates the fundamental group $\pi_1(S(2\gamma\oplus \R)) = \Z/2$, and $M_A\setminus \nu_0$ is the compact 4-manifold constructed in Section \ref{Section CappellShaneson}. 
Theorem \ref{Theorem CSF} implies that $Y$ is homeomorphic to $S(2\gamma\oplus \R)$. There are two $\Pin^+$-structures related to the $\eta$-invariant that we need to consider since $H^1(Y; \Z/2) = \Z/2$ \cite[Corollary 6.4]{[Stolz]}. We have that\begin{equation}\eta(S(2\gamma\oplus \R), \phi_i) = 0 \mod 2\Z\end{equation}for $i = 1, 2$ given that $S(2\gamma \oplus \R)$ is a $\Pin^+$-boundary of a $D^3$-bundle over $\R P^2$. Applying Theorem \ref{Theorem Stolz} we determine the value\begin{equation}\eta(Y, \phi_i^Y) = 1 \mod 2\Z\end{equation}for $i = 1, 2$. Since the values of the $\eta$-invariants are different for all $\Pin^+$-structures, we conclude that $Y$ is neither diffeomorphic nor stably diffeomorphic to $S(2\gamma\oplus \R)$ by Theorem \ref{Theorem Stolz Main}. The 2-sphere fiber $\Sigma \hookrightarrow S(2\gamma\oplus \R)$ has a tubular neighborhood $\nu (\Sigma) = D^2\times S^2$ and it is disjoint from the seam of the gluing along loops (\ref{Construction Y}) in the assembly of $Y$. 
We abuse notation and conclude that there is an embedded 2-sphere $\Sigma \subset Y$ with the desired properties.


\end{proof}

\subsection{Inequivalent smooth structures on non-orientable 4-manifolds.}\label{Section CSRP4}Cappell-Shaneson constructed a 4-manifold\begin{equation}\label{CS Manifold}Q:= \R P^4\#_{S^1} M_A = (D^2\widetilde{\times} \R P^2) \cup (M_A\setminus \nu_0)\end{equation}that is homeomorphic, but not diffeomorphic to $\R P^4$ \cite{[CappellShaneson], [Ruberman], [FreedmanQuinn]}. Let $\gamma\subset Q$ be a loop whose homotopy class generates the group $\pi_1Q$. A close look to Cappell-Shaneson's proof reveals their construction of a compact 4-manifold\begin{equation}\label{Exotic Compact Bundle}N: = Q\setminus \nu(\gamma)\end{equation} that is homeomorphic but not diffeomorphic to $D^2\widetilde{\times} \R P^2$ and whose boundary is $\partial N = S^2\widetilde{\times} S^1$ \cite[Section 2]{[Akbulut1]}. A handlebody of (\ref{Exotic Compact Bundle}) was drawn by Akbulut in \cite[Figures 1.26 and 1.33]{[Akbulut1]}. The second author of this note observed in \cite[Theorem 2]{[Torres]} that (\ref{Exotic Compact Bundle}) along with Theorems \ref{Theorem CSF} and \ref{Theorem Stolz Main} can be used to produce an inequivalent smooth structure on a myriad of closed non-orientable 4-manifolds, which include\begin{equation}\label{Smooth Structure Y}Y = S(2\gamma\oplus \R)\#_{S^1} M_A.\end{equation}and\begin{equation}\label{Inequivalent 1}Z = (\R P^4\#_{S^1} \R P^4)\#_{S^1} M_A = \R P^4\#_{S^1} Q.\end{equation}

\subsection{A result of Akbulut}\label{Section Result Akbulut}For non-orientable 4-manifolds with a $\Pin^+$-structure, the building block (\ref{Exotic Compact Bundle}) and the $\eta$-invariant are effective tools to produce inequivalent smooth structures (cf. \cite{[FintushelStern1], [FintushelStern2], [Kreck]}). The situation where there is no $\Pin^+$-structure available is strikingly different as the following theorem suggests.

\begin{theorem}\label{Theorem Akbulut}Akbulut \cite[Section 2]{[Akbulut1]}. Let $\alpha \subset \R P^2\times S^2$ be a simple loop whose homotopy class generates the group $\pi_1(\R P^2\times S^2)$. The 4-manifold\begin{equation}Y':= ((\R P^2\times S^2)\setminus \nu(\alpha))\cup (M_A\setminus \nu_0)\end{equation} is diffeomorphic to $\R P^2\times S^2$.
\end{theorem}

Theorem \ref{Theorem Akbulut} is one of the two key ingredients in the proof of Theorem \ref{Theorem A}.

\subsection{Connected sums with $\mathbb{CP}^2$}\label{Section CP2}The closed smooth 4-manifolds with fundamental group $\Z/2$ and Euler characteristic equal to two that are considered in this note compose three homeomorphism classes and five diffeomorphism classes. All these 4-manifolds become diffeomorphic under the connected sum of a single copy of the complex projective plane. The following proposition addresses the situation regarding the representatives of the three homeomorphism classes.

\begin{proposition}\label{Proposition CP2 Diffeomorphisms}There are diffeomorphisms\begin{equation}\label{Diffeomorphism 2}f_1: S(2\gamma \oplus \R) \# \mathbb{CP}^2\rightarrow (\R P^2\times S^2)\# \mathbb{CP}^2 \end{equation}and\begin{equation}\label{Diffeomorphism 3}f_2: (\R P^4\#_{S^1} \R P^4) \# \mathbb{CP}^2\rightarrow (\R P^2\times S^2)\# \mathbb{CP}^2.\end{equation}
\end{proposition}

\begin{proof}The existence of the diffeomorphism (\ref{Diffeomorphism 2}) follows from Lemma \ref{Lemma Gluck Twist} and the known property that Gluck twists come undone by taking a connected sum with a copy of $\mathbb{CP}^2$ \cite[Exercise 5.2.7 (b)]{[GompfStipsicz]}, \cite[Theorem 1.1]{[AkbulutYasui]}. The existence of the diffeomorphism (\ref{Diffeomorphism 3}) follows from Example \ref{Example Basic}, Lemma \ref{Lemma Stabilize Twisting RP2} and (\ref{Diffeomorphism 2}. 

\end{proof}

The following result implies that the inequivalent smooth structures that we have considered in this paper become diffeomorphic after adding a copy of $\mathbb{CP}^2$.

\begin{theorem}\label{Theorem CP2}Akbulut \cite[Theorem 2]{[Akbulut1]}. Let $Q$ be Cappell-Shaneson's exotic real projective 4-space that was constructed in (\ref{CS Manifold}). There is a diffeomorphism between $Q\#\mathbb{CP}^2$ and $\R P^4\# \mathbb{CP}^2$.
\end{theorem}

The existence of topologically equivalent and smoothly inequivalent 2- and 3-spheres in the non-trivial $\R P^2$-bundle over $S^2$ is an immediate consequence of Theorem \ref{Theorem CP2} and work of Cappell-Shaneson \cite{[CappellShaneson]}.

\begin{example}\label{Example Akbulut} Akbulut, Cappell-Shaneson. There is a pair $\{S^k_1, S^k_2\}$ of $k$-spheres for $k\in \{2, 3\}$ that are smoothly embedded in $\R P^4\# \mathbb{CP}^2$, whose complements have fundamental group of order two and such that there is a homeomorphism of pairs\begin{equation}(\R P^4\# \mathbb{CP}^2, S^k_1)\rightarrow (\R P^4\# \mathbb{CP}^2, S^k_2)\end{equation} but no such diffeomorphism of pairs exists. 
\end{example}

In order to prove Theorem \ref{Theorem CP2}, Akbulut shows that there is a diffeomorphism\begin{equation}\label{Diffeomorphism AKbulutN}N\#\mathbb{CP}^2 \rightarrow (D^2\widetilde{\times} \R P^2)\# \mathbb{CP}^2\end{equation} for the compact 4-manifold (\ref{Exotic Compact Bundle}). Corollary \ref{Corollary Akbulut} follows immediately from Akbulut's Theorem \ref{Theorem CP2} by noticing that both constructions (\ref{Smooth Structure Y}) and  (\ref{Inequivalent 1}) contain a copy of $N$. The latter is responsible for the inequivalent smooth structure. For example,\begin{equation}N\subset Z = ((\R P^4\#_{S^1}\R P^4)\#_{S^1} M_A) = \R P^4\#_{S^1}Q.\end{equation} We finish this section by stating a consequence of Theorem \ref{Theorem CP2} that is one of the main ingredients in the proof of Theorem \ref{Theorem C}: the inequivalent smooth structures (\ref{Smooth Structure Y}) and (\ref{Inequivalent 1}) are $\mathbb{CP}^2$-stable in the following sense.

\begin{corollary}\label{Corollary Akbulut}There are diffeomorphisms\begin{equation}f_3: Y \# \mathbb{CP}^2\rightarrow (\R P^2\times S^2)\# \mathbb{CP}^2\end{equation}and\begin{equation}f_4: Z\# \mathbb{CP}^2\rightarrow (\R P^2\times S^2)\# \mathbb{CP}^2.\end{equation}
\end{corollary}

\section{Proofs of main results}\label{Section Proofs Main Results}

\subsection{Knotted 2-spheres: Proof of Theorem \ref{Theorem A}}Construct\begin{equation}\label{Gluck Twist Y}Y':= (Y\setminus \nu (\Sigma)) \cup_{\varphi} (D^2\times S^2)\end{equation} by performing a Gluck twist to the 2-sphere $\Sigma \subset Y$ of Proposition \ref{Proposition Sphere Bundle}. Let $\alpha\subset \R P^2\times S^2$ be a simple loop that generates $\pi_1(\R P^2\times S^2) = \Z/2$. Another construction of (\ref{Gluck Twist Y}) is $Y'= ((\R P^2\times S^2)\setminus \nu(\alpha)) \cup (M_A\setminus \nu_0)$, and Theorem \ref{Theorem Akbulut} says that it is diffeomorphic to $\R P^2\times S^2$. From (\ref{Gluck Twist Y}), we see that there is a smoothly embedded 2-sphere $S \subset \R P^2\times S^2$ such that one recovers $Y$ by performing a Gluck twist on it. To conclude that the exteriors are not homotopy equivalent, we look at their fundamental groups. The exterior of the factor 2-sphere $(\R P^2\times S^2)\setminus \nu(\{pt\}\times S^2)$ is $Mb\times S^2$ and $(\R P^2\times S^2)\setminus \nu(S) = (Mb\times S^2)\#_{S^1} M_A$, where the circle sum is taken with matching $\Pin^+$-structures. These two compact 4-manifolds have non-isomorphic fundamental groups \cite[Proposition 2.4]{[CappellShaneson]}.

\hfill $\square$

\subsection{Surgeries to Cappell-Shaneson's mapping torus and smooth structures: Proof of Theorem \ref{Theorem B}}\label{Section Proof of Theorem B}Let $\alpha_0^2\subset M_A$ be the simple loop that represents the element $[\alpha_0^2]\in \pi_1(M_A)$ where $\alpha_0\subset M_A$ is the cross-section of the mapping torus (\ref{Mapping Torus}) that was described in Section \ref{Section CappellShaneson}. Carve out a tubular neighborhood $\nu(\alpha_0^2) = S^1\times D^3$ of the loop $\alpha_0^2\subset M_A$ and build the 4-manifold\begin{equation}\label{Construction CS Fr1}Y(\psi):= (M_A\setminus \nu(\alpha_0^2))\cup_{\psi} (D^2\times S^2),\end{equation} where $\psi = \id$ corresponds to the canonical framing defined in Section \ref{Section CappellShaneson}, and $\psi = \varphi$ corresponds to the non-canonical one, where the diffeomorphism $\varphi$ is described in (\ref{Diffeo Gluck}). Notice that the 4-manifold (\ref{Construction CS Fr1}) is a boundary component of the cobordism that is obtained from adding a 2-handle to $M_A\times I$. Rephrase (\ref{Construction CS Fr1}) as\begin{equation}\label{Construction CS FrG}Y(\psi) = (M_A\setminus \nu(\alpha_0^2))\cup_{\psi} (D^2\times S^2) = ((Mb\times S^2)\#_{S^1} M_A))\cup_{\psi} (D^2\times S^2).\end{equation}Lemma \ref{Lemma Gluck Twist} and Theorem \ref{Theorem Akbulut} imply that $Y(\id)$ is diffeomorphic to $\R P^2\times S^2$; notice that $Y(\id)$ is denoted by $Y'$ in Theorem \ref{Theorem Akbulut}. The 4-manifold $Y(\varphi)$ is denoted by $Y$ in Proposition \ref{Proposition Sphere Bundle} and Theorem \ref{Theorem A}.



\hfill $\square$

\begin{remark}\label{Remark Spheres} A Cappell-Shaneson homotopy 4-sphere is constructed by doing surgery along a cross-section of a certain 3-torus bundle over the circle $M_{A^2}$ with monodromy $A^2$ for $A\in \Gl(3; \Z)$ \cite{[CappellShaneson], [AitchisonRubinstein]}. A myriad of these homotopy 4-spheres are known to be diffeomorphic to $S^4$ by work of several authors \cite{[AitchisonRubinstein], [Akbulut4], [AkbulutKirby1],[AkbulutKirby2], [Gompf], [Gompf3], [Gompf2]}. Theorem \ref{Theorem B} addresses the corresponding situation in the non-orientable realm.

\end{remark}

\subsection{Surgery on $\R P^2$ and twisted doubles: Proof of Theorem \ref{Theorem R}}\label{Section Proof of Theorem R} Deconstruct the smooth structure on $S(2\gamma \oplus \R)$ of Proposition \ref{Proposition Sphere Bundle} as\begin{equation}\label{Decomposition Y}Y = N\cup_{\id} (D^2\widetilde{\times} \R P^2),\end{equation}where $N$ is the 4-manifold (\ref{Exotic Compact Bundle}).  A tubular neighborhood $\nu(\R P^2) = D^2\widetilde{\times} \R P^2$ of a smoothly embedded $\R P^2\subset Y$ is visible from (\ref{Decomposition Y}). Twisting this $\R P^2$ yields\begin{equation}\label{Decomposition Z}Z = N\cup_{\varphi'} (D^2\widetilde{\times} \R P^2). \end{equation}  

It is immediate that (\ref{Decomposition Z}) is homeomorphic to $\R P^4\#_{S^1}\R P^4$. Using Theorem \ref{Theorem Stolz Main} and the additivity of the $\eta$-invariant under the operation of circle sum, we conclude that $Z$ is neither diffeomorphic nor stably diffeomorphic to $\R P^4\#_{S^1}\R P^4$. 
\hfill $\square$







\subsection{Inequivalent smooth structures on non-orientable 4-manifolds with fundamental group $\Z/2\times \Z/2$: Proof of Theorem \ref{Theorem New Examples Structures}}\label{Section New Structures}The 4-manifold $M$ is homotopy equivalent to either $V = \R P^2\widetilde{\times}\R P^2$ or to the 4-manifold $W$ of Example \ref{Example New Construction} by \cite[Theorem 9.3]{[HambletonHillman]}. Follow (\ref{Decomposition Bundle}) and (\ref{Decomposition Not Bundle}) to construct\begin{center}$V' = N\cup_{\id} (Mb\widetilde{\times} \R P^2)$ and $W' = N\cup_{\varphi'} (Mb\widetilde{\times} \R P^2)$,\end{center}where $N$ is the 4-manifold (\ref{Exotic Compact Bundle}). There are homeomorphisms $V'\rightarrow V$ and $W'\rightarrow W$ by construction. To see that these 4-manifolds are not diffeomorphic, we show that they are not $\Pin^+$-cobordant by computing the $\eta$-invariant for all choices of $\Pin^+$-structure \cite{[Stolz]} (cf. Theorem \ref{Theorem Stolz Main} and Remark \ref{Remark Theorem Stolz Main}). Since $H^1(M; \Z/2) = \Z/2\times \Z/2$, there are four $\Pin^+$-structures to consider. We see that\begin{center}$\eta(V, \phi^V_i) = \eta(W, \phi^W_i) = \pm \frac{1}{8} \mod 2\Z$\end{center} for $i = 1, 2, 3, 4$. These values are computed as follows. Let $\alpha\subset \R P^4$ be a simple loop whose homotopy class generates the group $\pi_1(\R P^4) = \Z/2$ and let $\gamma\subset \R P^3\widetilde{\times} S^1$ be a cross-section of the non-trivial $\R P^3$-bundle over $S^1$. Glue $\R P^4$ and $\R P^3\widetilde{\times} S^1$ along these loops and with matching $\Pin^+$-structures to produce\begin{center}$V = (D^2\widetilde{\times} \R P^2)\cup_{\id} (\R P^3\widetilde{\times} S^1\setminus \nu(\gamma))$ and $W = (D^2\widetilde{\times} \R P^2)\cup_{\varphi'} (\R P^3\widetilde{\times} S^1\setminus \nu(\gamma))$\end{center} with all their possible $\Pin^+$-structures \cite{[Kirby], [KirbyTaylor]}. Lemma \ref{cobordism} says that there are $\Pin^+$-cobordisms between $(V, \phi^V_i)$ or $(W, \phi^W_i)$ and the disjoint union $(\R P^4, \pm \phi^{\R P^4})\sqcup (\R P^3\widetilde{\times} S^1, \phi^{\R P^3\widetilde{\times} S^1}_i)$. The additivity of the $\eta$-invariant with respect to disjoint unions implies that the values hold. Analogously, we get $\eta(V', \phi^{V'}_i) = \eta(W', \phi^{W'}_i) = \pm \frac{7}{8} \mod 2\Z$ for $i = 1, 2, 3, 4$ by employing $(Q, \pm \phi^Q)$ \cite{[Stolz]}. 
\hfill $\square$

\begin{remark}\label{Remark New Twisting}Follow the assembly of $N$ in (\ref{CS Manifold}) to construct\begin{center}$V' = (D^2\widetilde{\times} \R P^2)\cup_{\id} ((Mb\widetilde{\times} \R P^2)\cup (M_A\setminus \nu_0))\setminus \nu(\alpha)$\end{center} and \begin{center}$W' = (D^2\widetilde{\times} \R P^2)\cup_{\varphi'} ((Mb\widetilde{\times} \R P^2)\cup M_A\setminus \nu_0))\setminus \nu(\alpha)$,\end{center}where $\alpha\subset (Mb\widetilde{\times} \R P^2)\cup (M_A\setminus \nu_0)$ is a suitably chosen simple loop. Therefore, $V'$ is obtained by twisting a smoothly embedded $\R P^2\subset W'$ and vice versa.

\end{remark}

\subsection{Knotted 2-spheres: Proof of Theorem \ref{Theorem C}}We select the five 2-spheres $\{S_i: i = 0, \ldots, 4\}$ as follows. Denote by $S_0$ the image of the canonical embedding of $\mathbb{CP}^1\hookrightarrow (\R P^2\times S^2)\# \mathbb{CP}^2$ for $\mathbb{CP}^1\subset \mathbb{CP}^2\setminus D^4$. Set $S_i:= f_i(\mathbb{CP}^1)$ where $f_i$, depending on $i \in \{1, 2, 3, 4\}$, is the corresponding diffeomorphism of Proposition \ref{Proposition CP2 Diffeomorphisms} and Corollary \ref{Corollary Akbulut}. Item (1) and Item (3) of Theorem \ref{Theorem C} follow by construction. Item (2) follows from the fact that, by examining the diffeomorphisms $f_1$ and $f_2$ at the level of handlebodies, one can explicitly see that the homology classes $[S_1]$ and $[S_2]$ are both equal to $[\mathbb{CP}^1] + [\{pt\} \times S^2] \in H_2((\R P^2 \times S^2)\# \mathbb{CP}^2; \Z)$. Indeed, a handlebody for $S(2 \gamma \oplus \mathbb{R})\# \mathbb{CP}^2$ is obtained by adding a 1-framed unknot to the handlebody of $M_1 = S(2\gamma \oplus \mathbb{R})$ given on top of Figure 1 for $k = 1$. The 2-sphere $\mathbb{CP}^1 \subset \mathbb{CP}^2\setminus D^4\subset S(2\gamma \oplus \R)\#\mathbb{CP}^2$ arises as the union of the core of the 1-framed 2-handle and a 2-disk in $S^3$ bounding its attaching circle. In the following, we keep track of the image of this 2-sphere via the diffeomorphism $f_1$. Slide the (1,0)-framed 2-handle over the 1-framed 2-handle. Slide the 1-framed 2-handle over the 0-framed 2-handle to unlink it from the diagram. These two handleslides yield a handlebody depicted on top of Figure 1 with $k = 2$ and a 1-framed unknot, which is a handlebody for $(\R P^2\times S^2)\#\mathbb{CP}^2$. In particular, the 2-sphere $\mathbb{CP}^1 \subset S(2\gamma \oplus \R) \# \mathbb{CP}^2$ is mapped onto a 2-sphere that is obtained by tubing a factor 2-sphere in $\R P^2\times S^2$ with $\mathbb{CP}^1 \subset (\mathbb{R}P^2 \times S^2)\# \mathbb{CP}^2$. This shows that $[S_1]=[S_0]+[\{pt\} \times S^2] \in H_2((\R P^2\times S^2)\#\mathbb{CP}^2; \Z)$. The claim about the homology class of $S_2$ is proven in a similar way. On the other hand, one can check that the map induced in homology by the diffeomorphism $N \# \mathbb{CP}^2 \approx (D^2 \widetilde \times \mathbb{R}P^2)\# \mathbb{CP}^2$ described in (\ref{Diffeomorphism AKbulutN}) sends the class of the exceptional sphere $\mathbb{CP}^1 \subset \mathbb{CP}^2$ onto itself. The decompositions (\ref{S_3}) and (\ref{S_4}) hence imply the statement for the homology classes $[S_3], [S_4] \in H_2((\R P^2\times S^2)\#\mathbb{CP}^2; \Z)$, by reducing to the case of $[S_1]$ and $[S_2]$ respectively. Item (4) is immediate since\begin{center}$(\R P^2\times S^2)\setminus D^4$, $S(2\gamma\oplus \R)\setminus D^4$ and $(\R P^4\#_{S^1}\R P^4)\setminus D^4$\end{center} are pairwise non-homotopy equivalent. Item (5) follows from Corollary \ref{Corollary Akbulut} and Proposition \ref{Proposition Sphere Bundle}. Indeed, a diffeomorphism of pairs\begin{center}$((\mathbb{R}P^2 \times S^2)\# \mathbb{CP}^2, S_1) \to ((\mathbb{R}P^2 \times S^2)\# \mathbb{CP}^2, S_3)$\end{center} would induce a diffeomorphism between the complements of the 2-spheres\begin{center}$(\mathbb{R}P^2 \times S^2)\# \mathbb{CP}^2 \setminus \nu(S_1) \to (\mathbb{R}P^2 \times S^2)\# \mathbb{CP}^2 \setminus \nu(S_3)$.\end{center} This is not possible: Theorem \ref{Theorem A} says that the complements $(\mathbb{R}P^2 \times S^2)\# \mathbb{CP}^2 \setminus \nu(S_1) \approx S(2\gamma \oplus \R) \setminus D^4$ and $(\mathbb{R}P^2 \times S^2) \# \mathbb{CP}^2 \setminus \nu(S_3)\approx Y \setminus D^4$ are not diffeomorphic. In a similar way, the existence of a diffeomorphism of pairs $((\mathbb{R}P^2 \times S^2)\# \mathbb{CP}^2, S_2) \to ((\mathbb{R}P^2 \times S^2)\# \mathbb{CP}^2, S_4)$ would contradict Theorem \ref{Theorem R}.

\hfill $\square$

\begin{figure}{\label{Figure 1}}
\begin{center}
\includegraphics[width=100mm,scale=0.4]{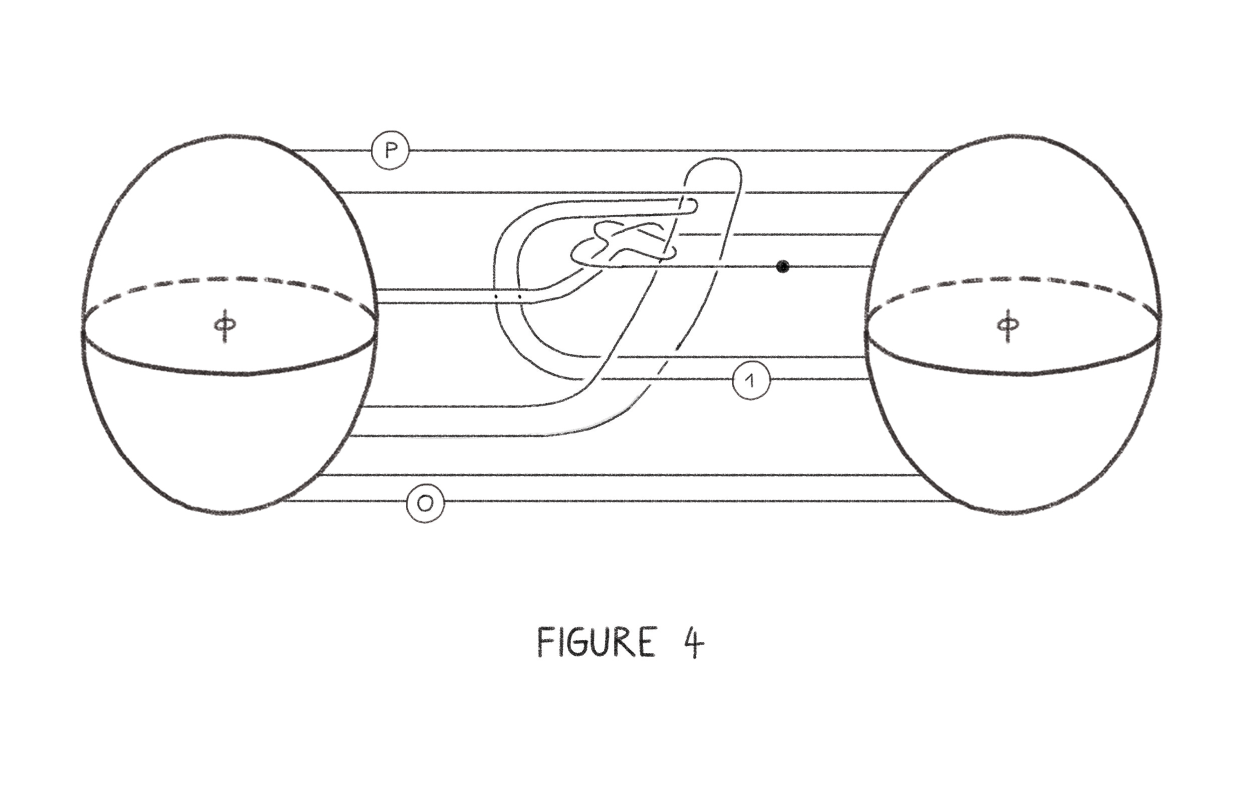}
\captionsetup{labelformat=empty}
\caption{Handlebodies of (\ref{Decomposition Y}) and (\ref{Decomposition Z}) are drawn in Figure 4. The value $p = - 1$ corresponds to $Y$, and $p = 1$ to $Z$. Figure 4 corrects \cite[Fig 2]{[Torres]}.}
\end{center}
\end{figure}



\end{document}